\documentclass{ifacconf}

\usepackage{amsmath}
\usepackage{amsfonts}
\usepackage{amssymb}
\usepackage{graphicx}
\usepackage[usenames]{color}
\usepackage{natbib}   
\usepackage{flushend}
\newtheorem{theorem}{Theorem}
\newtheorem{defin}{Definition}

\newtheorem{remark}{Remark}
\newcommand{\pder}[2][]{\frac{\partial#1}{\partial#2}}

\newcommand{\RR}{\mathbb{R}}

\begin{document}

\begin{frontmatter}

\title{Control Contraction Metrics \\and Universal Stabilizability} % Title, preferably not more than 10 words.

\thanks[footnoteinfo]{This work was supported in part by the Australian Research Council.}

\author[First]{Ian R. Manchester} 
\author[Second]{Jean-Jacques E. Slotine} 

\address[First]{ACFR, School of Aerospace, Mechanical and Mechatronic Engineering, University of Sydney, Australia.\\ ian.manchester@sydney.edu.au}                                              
\address[Second]{Nonlinear Systems Laboratory, Dept. of Mechanical Engineering,\\ Massachusetts Institute of Technology}

%\begin{keyword}                           % Five to ten keywords,  
%Cicero; Catiline; orations.               % chosen from the IFAC 
%\end{keyword}                             % keyword list or with the 
                                          % help of the Automatica 
                                          % keyword wizard

\begin{abstract}    
In this paper we introduce the concept of universal stabilizability: the condition that every solution of a nonlinear system can be globally stabilized. We give sufficient conditions in terms of the existence of a control contraction metric, which can be found by solving a pointwise linear matrix inequality. Extensions to approximate optimal control are straightforward. The conditions we give are necessary and sufficient for linear systems and certain classes of nonlinear systems, and have interesting connections to the theory of control Lyapunov functions.
\end{abstract}

\end{frontmatter}

%%%%%%%%%%%%%%%%%%%%%%%%%%%%%%%%%%%%%%%%%%%%%%%%%%%%%%%%%%%%%%%%%%%%%%%%%%%%%%%%
\section{Introduction}

Constructive control design for nonlinear systems remains a challenging problem, even in the case of full state feedback \citep{slotine1991applied},  \citep{isidori1995nonlinear}, \citep{kokotovic2001constructive}.

The classical Lyapunov stability theory leads to necessary and sufficient conditions in terms of existence of control Lyapunov functions, however these may be difficult to find. A full state feedback linearization allows linear control tools to be applied; however finding such a transform can be non-trivial even when it can be shown that one exists \citep{isidori1995nonlinear}.

Constructive design tools tend to be more limited in their scope of application. Backstepping and related methods provide a systematic and constructive procedure \citep{krstic1995nonlinear}, but are generally limited to systems of a particular ``triangular'' structure. For certain mechanical and electrical systems, energy and passivity methods can be used \citep{schaft1999l2}.

Recently, there has been increased interest in methods for systems analysis and control design based on convex optimization, for example linear matrix inequalities, integral quadratic constraints, and sum-of-squares programming \citep{boyd1994linear}, \citep{megretski1997system}, \citep{parrilo2000structured}. For nonlinear control design, the density functions of \citep{rantzer2001dual}, \citep{prajna2004nonlinear} and related techniques of occupation measures \citep{lasserre2008nonlinear} and control Lyapunov measures \citep{vaidya2010nonlinear} explicitly address convexity of criteria. Another approach is to piece together stabilized trajectories, with regions of stability verified via sum-of-squares programming \citep{tedrake2010lqr}.

Constructing a control-Lyapunov or density function requires prior knowledge of the solution to be stabilized. However, in many applications with complex system architectures, the role of feedback control is to track a setpoint or trajectory that changes in real-time. For nonlinear systems the problem of stabilizing changing trajectories can be quite different to stabilizing a single {\em a priori} known trajectory.

An alternative to explicit design is model predictive control, based on real-time numerical optimization. For many years this method was limited to relatively slow linear processes, but rapid advances in computer technology and optimization algorithms mean it is emerging as a feasible tool for nonlinear problems (see, e.g. \cite{diehl2009efficient}). Despite some clear benefits, it generally remains difficult to predict or analyse performance of nonlinear MPC schemes by any method other than exhaustive simulations.

In this paper, we introduce the concept of {\em universal stabilizability}, i.e. the property that {\em every} solution of a system is globally stabilizable. Using the tools of contraction analysis \citep{Lohmiller98}, we give convex sufficient conditions for universal stabilizability, and a constructive control design procedure. Useful properties such as exponential stability and approximate linear quadratic optimality are natural extensions, and for the case of linear systems, our conditions reduce to well-known necessary and sufficient conditions.

Historically, basic convergence results on contracting systems can be traced back to the results of \cite{lewis1949metric} in terms of Finsler metrics, and results of  \cite{hartman1961stability} and  \cite{demidovich1962dissipativity}. Extensions to analysis of limit cycles have recently been developed by the authors \citep{manchester2013transverse} building upon the early work of \cite{Borg} and \cite{hartman1962global}. Connections with Finsler structures and Lyapunov theory have recently been explored in detail by \cite{forni2012differential}.

The conditions we give are state-dependent linear matrix inequalities. The conditions are formally similar to those studied for a class of nonlinear $H^\infty$ problems in \cite{lu1995} and to state-dependent Riccati equations (see, e.g., \cite{cloutier1997state}), however the class of systems and the explicit control construction we consider are quite different. An interesting feature of our method is that it breaks the control-design problem into two stages: computation of a metric, which can be performed in advance of knowledge of the solution to be stabilized, and an on-line path integration to compute the control signal. 

Our methods build significantly upon the control design suggested in \cite{lohmiller1998contraction}. In that paper, it was argued that requiring a nonlinear system to be feedback transformable to a stable linear system -- as in feedback linearization --
was too strong, since it requires an involutivity condition to be added to the rather natural controllability
condition \citep{isidori1995nonlinear}. If instead one required only that the system be feedback transformable to a ``nice'' stable nonlinear system
(a contracting system), then this can be achieved using the controllability condition alone. 

\section{Universal Stabilizability}

For most of this paper, we will consider a nonlinear time-dependant control-affine system
\begin{equation}\label{eq:sys}
\dot x(t) = f(x(t),t) +B(t)u(t)
\end{equation}
where $x(t)\in\RR^n, u(t)\in\RR^m$ are state and control, respectively, at time $t\in\RR^+:=[0,\infty)$. The function $f: \RR^n\times \RR^+ \rightarrow \RR^n$ is assumed to be smooth, and $B:\RR^+ \rightarrow \RR^{n\times m}$ is a time-dependent matrix.

Contraction analysis is the study of \eqref{eq:sys} by way of the associated system of differential dynamics:
\begin{equation}\label{eq:diffdyn}
\dot\delta_x(t) = A(x,t)\delta_x(t)+B(t)\delta_u(t)
\end{equation}
where $A(x,t) = \pder{x}f(x,t)$ is the Jacobian matrix.

In Section \ref{sec:gen} we will consider more general systems $\dot x= f(x,u,t)$ but for the moment we note that many systems not naturally appearing in the form \eqref{eq:sys} can be put in that form, either exactly or approximately, by change of variables or introducing new states.

A solution of \eqref{eq:sys} is a pair of vector signals $(x(t), u(t))$ satisfying \eqref{eq:sys} over the interval $\RR^+$. For simplicity, in this paper we will assume that \eqref{eq:sys} is such that solutions exist and are unique. As such, we will also use the notation $\phi(t, x_0, u)$ to denote the solution of $x$ at time $t>0$ of \eqref{eq:sys} starting from initial condition $x(0) = x_0$ and under the application of the control input $u(\tau), \tau \in [0, t]$.

A static state-feedback controller is a function $k:\RR^n \rightarrow \RR^m$. The main objective is to design such a function so that the behaviour of the closed-loop system
\begin{equation}\label{eq:clsys}
\dot x(t) = f(x(t),t)+B(t)k(x(t),t)
\end{equation}
is in some sense desirable, e.g.  globally stable or optimal in some sense. Solutions of the closed-loop system with a given state-feedback controller $k$ will be denoted $\phi_k(t, x_0)$.

As is standard, a solution $(x^\star, u^\star)$ defined on $[0, \infty)$ is said to be globally asymptotically stabilized by a feedback controller $u=k(x,t)$ if both of the following hold:
\begin{enumerate}
\item For any $\alpha$ there exists an $\epsilon$ such that $|x_0-x^\star(0)|<\epsilon$ implies $|\phi_k(t, x_0)-x^\star(t)|<\alpha$,
\item For any initial condition $x_0\in \RR^n$, the closed loop solution satisfies $|\phi_k(t, x_0)-x^\star(t)|\rightarrow 0$.
\end{enumerate}
Global exponential stabilization refers to the stronger condition that there exists a $K$ and $\lambda$ such that \[|\phi_k(t, x_0)-x^\star(t)|\le K e^{-\lambda t}|x_0-x^\star(0)|\] for all $x(0)$.

In this paper, we will study the following property:
\begin{defin} A system of the form \eqref{eq:sys} is said to be {\em universally stabilizable by state feedback} if for any solution $(x^\star, u^\star)$ defined on $t\in [0, \infty)$ there exists a state feedback controller $k:\RR^n\rightarrow \RR^m$ such that $(x^\star, u^\star)$ is globally stabilized by $u=k(x,t)$. 
\end{defin}
Analogously, we also consider the notion of {\em universally exponentially stabilizable with rate $\lambda$}.

Note that universal stabilizability is a stronger condition than global stabilizability of a particular solution.

\section{Control Contraction Metrics}

The following definition is central to this paper:

\begin{defin}A function $V(x, \delta_x, t)=\delta_x'M(x,t)\delta_x$, with $c_1I\le M(x,t) \le c_2I$ for some $c_2\ge c_1 >0$, is said to be a {\em control contraction metric} for the system \eqref{eq:sys} if $\pder[V]{x}B(t)=0$ and 
\begin{equation}\label{eq:ctrbmet}
\pder[V]{\delta_x}B(t)=0 \Longrightarrow \pder[V]{t} + \pder[V]{x}f(x,t)+\pder[V]{\delta_x}A(x,t)\delta_x<0
\end{equation}
for all $x, t$.
\end{defin}

We will also make use of a Riemannian distance function between any two points at a given time $d(x_1, x_2, t): \RR^n \times \RR^n \times \RR^+ \rightarrow \RR^+$ defined like so: let $\Gamma(x_1, x_2)$ denote the set of all smooth paths connecting $x_1$ and $x_2$, where each $\gamma\in\Gamma(x_1,x_2)$ is parametrised by $s\in [0, 1]$, i.e. $\gamma(s):[0, 1] \rightarrow \RR^n$. The path length of $\gamma$ is then defined as
\[
L(\gamma,t) := \int_0^1 D\left(\gamma(s),\frac{\partial}{\partial s}\gamma(s),t\right)ds
\]
where $D$ is a Finsler function \citep{bao2000introduction}. In this paper, we will primarily choose $D(x, \delta_x, t) = \sqrt{V(x, \delta_x, t)}$ or $D(x, \delta_x, t)=|\delta_x|$. The distance between two points is then defined as
\[
d(x_1, x_2, t) = \min_{\gamma\in\Gamma(x_1, x_2)}L(\gamma,t)
\]
The existence of a minimizing path, which we denote $\gamma_{x_1}^{x_2}(t,s)$,  is implied by the Hopf-Rinow Theorem.

\begin{remark}
In this conference paper we restrict ourselves to simple choices of $V$ and $D$. Other choices are possible, including non-quadratic metrics and matrix measures. See, e.g., \cite{lohmiller1998contraction} and \cite{forni2012differential} for detailed discussion of non-quadratic metrics in contraction analysis.
\end{remark}

The main utility of a control contraction metric is that it allows one to construct stabilizing controls by computing path integrals of a suitable $\delta_u$.

\begin{theorem} If a control contraction metric exists for a system of the form \eqref{eq:sys}, then the system is universally stabilizable by static state feedback.
\end{theorem}
\textbf{Proof:} 
Consider an arbitrary initial condition $x(0)\in\RR^n$ and a desired feasible trajectory $(x^\star(t), u^\star(t))$ defined on $t\in[0, \infty)$.

For each time $t$, consider the distance $d(x^\star(t), x(t), t)$. By the construction of the distance given above, if $d\rightarrow 0$ then $x(t)\rightarrow x^\star(t)$. So the objective of this proof is to show the existence of a control law such that $d\rightarrow 0$ as $t\rightarrow \infty$.

Consider also the minimal path $\gamma(t,s)=\gamma_{x_1}^{x_2}(t,s)$. For a given time $t$, at each point $s\in[0,1]$ along the path, consider two possibilities: firstly, suppose
\[
\pder{\delta_x}V(x, \delta_x, t)B(t) = 0
\]
when the left hand side is evaluated at $x = \gamma(t,s), \delta_x = \pder[\gamma(s,t)]{s}$. Then by the definition of a control contraction metric, one has
\[
\frac{d}{dt}V = \pder[V]{t} + \pder[V]{x}f(x,t)+\pder[V]{\delta_x}A(x,t)\delta_x<0.
\]
Secondly, suppose that 
\[
\pder{\delta_x}V(x, \delta_x, t)B(t) \ne 0
\]
at $x = \gamma(t,s), \delta_x = \pder[\gamma(s,t)]{s}$. Then 
\[
\frac{d}{dt}V = \pder[V]{t} + \pder[V]{x}f(x,t)+\pder[V]{\delta_x}(A(x,t)\delta_x+B(t)\delta_u)
\]
is affine in $\delta_u$, and therefore there exists $\delta_u$ making $\frac{d}{dt}V<0$. Choose any such $\delta_u$ and denote it $\delta_u^\star(t,s)$.

Let us construct a control signal $u(t,s)$ for each point on the path $\gamma(t,s)$ like so:
\[
\bar u(t,s) = u^\star(t)+\int_0^s \delta_u^\star(t, \mathfrak s)d\mathfrak s
\]
then we have $u(t, 0) = u^\star(t)$ and $\pder[\bar u(t,s)]{s} = \delta_u^\star(t,s)$

This implies that for all $(t, s)$ we have \[\frac{d}{dt}V\left(\gamma(s),\frac{\partial}{\partial s}\gamma(s),t\right)<0\] and therefore as $t\rightarrow \infty$, we have
$V\left(\gamma(s),\frac{\partial}{\partial s}\gamma(s),t\right)\rightarrow 0$ for all $s$, implying $D\left(\gamma(s),\frac{\partial}{\partial s}\gamma(s),t\right)\rightarrow 0$ for all $s$, further implying $d(x^\star(t), x(t))\rightarrow 0$. 

The specific control signal that is applied to the system is $\bar u(t,1)$. Note that this procedure defines a static state-feedback function $u=k(x,t)$.

$\Box$

\section{Convex Condition and Construction of Explicit Control}\label{sec:control}

The second main result of this paper is the following, where 
we use the notation $S_+^n$ to denote the convex cone of symmetric positive semidefinite matrices:
\begin{theorem}\label{Wthm}
Consider the system \eqref{eq:sys} with differential dynamics \eqref{eq:diffdyn}. If there exists a matrix function $W(x,t):\RR^n \times \RR \rightarrow S_+^n$, a function $\rho(x,t)\ge 0$, and constants $\alpha_2\ge \alpha_1>0$ such that
\begin{align}
\alpha_1 I\le W\le \alpha_2 I, \label{eq:Wcond}\\
-\dot W + WA' + AW -\rho BB'< 0,\label{eq:Wdotcond}
\end{align}
for all $x, t$, then the system \eqref{eq:sys} is universally stabilizable by static state feedback, and $\delta_xW(x,t)^{-1}\delta_x$ is a control contraction metric.
\end{theorem}

In the above condition, $\dot W(x,t)$ is a matrix with the $i, j$ element given by $\pder[W_{i,j}]{t}+\pder[W_{i,j}]{x}(f(x,t)+B(t)u)$.

\textbf{Proof:} 
We will show that by taking $M(x,t) = W(x,t)^{-1}$, the conditions of the theorem prove that $V(x, \delta_x,t)=\delta_x'M(x,t)\delta_x$ is a control contraction metric.

Indeed, make the substitution of $M$ for $W$ and multiply \eqref{eq:Wdotcond} on either side by $M$, the conditions of the theorem are equivalent to 
\begin{align}\label{eq:MdotCond}
\dot M + A'M + MA -\rho MBB'M&<0,
\end{align}
$\forall x, u$ with $\frac{1}{\alpha_2}I \le M(x,t)\le \frac{1}{\alpha_1}I$. Note that $\dot M(x,t) = -M(x,t)\dot W(x,t) M(x,t)$ follows from the total derivative of the identity $M(x,t)W(x,t) = I$ with respect to time.

%Now, at any time $t\ge 0$ and for any particular $x(t)\in\RR^n$, a smooth path connecting $x(t)$ to $x^\star(t)$ can be represented via a mapping $\gamma:[0,1]\rightarrow
%\RR^n$ with $\gamma(0) = x^\star(t)$ and $\gamma(1) = x(t)$,  parametrized so that $\frac{\partial \gamma(s)}{\partial s} \ne 0$ for all $s$.

%For each $t$, by $\Gamma_t(x^\star(t), x(t))$ the set of all such smooth paths connecting
%$x(t)$ to $x^\star(t)$ and associate with each $\gamma_t\in\Gamma_t(x^\star(t), x(t))$ a length
%\[
%L(\gamma) = \int_0^1 V\left(\gamma(s),\frac{\partial}{\partial s}\gamma(s),t\right)ds
%\]
%where 
%\[
%V(x,\delta,t) = \sqrt{\delta_x'M(x,t)\delta_x}.
%\]

Let $\gamma_t(s)$ be the minimizing path connecting $x^\star(t)$ to $x(t)$, parameterized by $s\in[0, 1]$ and consider the control signal generated by
\[
u(t) = u^\star(t) -\frac{1}{2}\int_0^1 \rho(\gamma_t(s),t)B(t)'M(\gamma_t(s),t)\pder[\gamma_t]{s} ds
\]
is stabilizing to $x^\star(t)$. 

Now, along the path $\gamma$, we have $\delta_u=-\frac{1}{2}B'M(x)\delta_x$ so
\begin{align}
\frac{d}{dt}\delta_x'M\delta_x =& \delta_x'\dot M\delta_x +2\delta'M[A\delta+B\delta_u]\notag\\
=&  \delta_x'[\dot M + A'M + MA-\rho MBB'M]\delta_x<0\label{eq:metric_decrease}
\end{align}
where the last inequality is implied by  \eqref{eq:MdotCond}.

%Now, since $\gamma$ is assumed to be minimizing, 
%\[
%\frac{d}{dt}d(x^\star(t),x(t)) \le \int_0^1 \left[\frac{d}{dt}V\left(\gamma_t(s),\pder[\gamma_t]{s}(s) \right)\right]ds.
%\]
%and the negativity of the integrand on the right-hand-side is implied by \eqref{eq:metric_decrease}, since $\dot V(x, \delta_x,t) = \frac{1}{2V(x, \delta_x,t)}\frac{d}{dt}(\delta_x'M(x,t)\delta_x)$ and $V$ is positive-definite.

%Condition \eqref{eq:Wcond} ensures that convergence in the Riemannian metric implies convergence in Euclidean metric $|x-x^\star|$.

$\Box$

\begin{remark} We note that the conditions are pointwise linear matrix inequalities in $W(x,t)$ and $\rho(x,t)$, and are thus amenable to solution via convex optimization algorithms. We will address this further in Section \ref{sec:comp}.
\end{remark}

For systems of the form \eqref{eq:sys}, the matrix $\dot W$ has elements $\dot W_{i,j}=\pder[W_{i,j}]{t}+\pder[W_{i,j}]{x}(f(x,t)+B(t)u)$, which are affine in the control signal $u$. Therefore for inequality \eqref{eq:Wdotcond} to hold for all $u$, it is necessary that $\pder[W_{i,j}]{x}B(t)=0$, i.e. $W(x,t)$ must be constant along the subspace $S_u$ of state space in which the control signal directly acts.

This property deserves further investigation, but we remark that if $f(x,t)$ is globally Lipschitz, then tracking errors in $S_u$ will generally be stabilized by sufficiently high-gain linear feedback as verified by a Lyapunov function that is quadratic with respect to $S_u$. This is closely connected to the method of sliding-mode control \citep{slotine1991applied}. Alternatively, dynamics in these directions can be directly cancelled, in a simple application of feedback linearization \citep{isidori1995nonlinear}.

We also note that the restriction $\pder[W_{i,j}]{x}B(t)=0$ is removed in Section \ref{sec:gen}, at the expense of a more complex control computation.

\subsection{Conditions for Exponential Stabilization}\label{sec:exp}

A stronger result is the following:
\begin{theorem}
Consider the system \eqref{eq:sys} with differential dynamics \eqref{eq:diffdyn}. If there exists a matrix function $W(x,t)\in S_+^n$, $\rho(x,t)\ge 0$ and  $\alpha_2\ge\alpha_1>0$ such that
\begin{align}
\alpha_1I\le W\le \alpha_2 I, \\
-\dot W + WA' + AW -\rho BB' &\le -2\lambda W(x),\label{eq:Wdotcondexp}
\end{align}
for all $x, u, t$, then the system is universally exponentially stabilizable with rate $\lambda$ by state feedback.
\end{theorem}

The proof uses the same controller construction and follows the same outline as the first theorem, so we provide a sketch here: by taking $M(x,t) = W(x,t)^{-1}$ the conditions of the theorem are equivalent to $0<M(x)<\frac{1}{\alpha}I$ and
\begin{align}
\dot M + A'M + MA& -MBB'M\le-2\lambda M.
\end{align}
and, under the action of the control,
\[
\frac{d}{dt}\left(\delta_x'M(x,t)\delta_x\right) \le -2\lambda \delta_x'M(x,t)\delta_x.
\]
So there exists a $K_1$ such that $V(s,t)\le K_1e^{-2\lambda t}V(s,0)$ for all $s$.
Now, from the definition of $V$ and $D$ we can give $c_2\ge c_1>0$ such that
$c_1 \sqrt{V(t,s)} \le D(t,s)\le c_2 \sqrt{V(t,s)}$ it follows that $D(t,s) \le K e^{-\lambda t}$ with $K=\frac{c_2}{c_1}\sqrt{K_1}$.

This further implies that \[ d(x^\star(t), x(t))\le Ke^{-\lambda t}d(x^\star(0), x(0))\] which completes the proof.

$\Box$

\section{Guaranteed Quadratic Cost Control}\label{sec:gcc}

Explicit solutions of nonlinear optimal control problems are generally extremely challenging to find. A general approach is to formulate solving a nonlinear partial differential equation, the Hamilton Jacobi Bellman equation, and attempt to find a so-called viscosity solution. Unfortunately, in most cases of practical interest this computational problem is intractible.

A less demanding approach, frequently used in robust control, is to find a {\em guaranteed cost control}  \citep{petersen2000robust}. In this framework, one tries to find a controller that guarantees a conservative upper bound on the cost, which is usually obtained via dissipation inequalities. Then one can approach approximate optimal control by minimizing this upper bound. In this section we show that our method can be extended to guaranteed cost control of nonlinear systems, by making use of {\em differential} dissipation inequalities \citep{manchester2013transverse}.

\begin{theorem}\label{QGCCthm}
Given $Q(t)>0$ and $R(t)>0$ and the differential dynamics \eqref{eq:diffdyn}, suppose there exists a $W(x,t)$ such that
\begin{equation}\label{eq:Wcondgcc}
\begin{bmatrix}
(\dot W -WA' -AW+BR^{-1}B') & W\\W & Q^{-1}\end{bmatrix} \ge 0
\end{equation}
for all $x, u, t$. Now, for any feasible trajectory $(x^\star(t), u^\star(t))$, consider the LQ cost function
\[
J=\int_0^\infty \begin{bmatrix} (x(t)-x^\star(t)) \\  (u(t)-u^\star(t)) \end{bmatrix} '\begin{bmatrix}Q(t) & 0\\0 & R(t) \end{bmatrix}\begin{bmatrix} (x(t)-x^\star(t)) \\  (u(t)-u^\star(t)) \end{bmatrix}  dt.
\]
then there exists a state-feedback controller such that
\[
J(x(0))\le \int_0^1 \delta_x'W(x)^{-1}\delta_x ds
\]
with $x=\gamma(t,s)$ and $\delta_x = \pder[\gamma(t,s)]{s}$ over the path $\gamma_{x^\star(0)}^{x(0)}(t,s)$.
\end{theorem}

\textbf{Proof:} Follows from similar argument as above and the following differential dissipation inequality
\[
\frac{d}{dt} (\delta_x M\delta_x)\le - \begin{bmatrix} \delta_x \\ \delta_u \end{bmatrix} '\begin{bmatrix}Q(t) & 0 \\0 & R(t) \end{bmatrix}\begin{bmatrix} \delta_x \\ \delta_u \end{bmatrix}
\]
along any path $\gamma$ joining $x^\star(t)$ and $x(t)$.

This inequality is satisfied by design of controllers satisfying the differential Riccati inequality
\[
\dot M + A'M + MA -MBR^{-1}B'M+Q\le 0
\]
which can be transformed to
\[
-\dot W + WA' + AW -BR^{-1}B' +WQW \le 0
\]
which can be linearized via Schur complement to
\[
\begin{bmatrix}
(\dot W -WA' -AW+BR^{-1}B') & W\\W & Q^{-1}\end{bmatrix} \ge 0.
\]

To show that  the differential dissipation inequality implies the cost bound, it is more convenient to consider the equivalent formulation
\[
J=\int_0^\infty  \left | H(t) \begin{bmatrix} (x(t)-x^\star(t)) \\  (u(t)-u^\star(t)) \end{bmatrix}\right|^2dt
\]
where $H(t)$ satisfies
\[H(t)'H(t) =  \begin{bmatrix}Q(t) & 0 \\0 & R(t) \end{bmatrix}\]
e.g. by Cholesky factorization. 

The upper bound for the cost can be obtained by noting that the Cauchy-Schwarz inequality on $L^2[0, 1]$ implies
\[
\int_0^1\left| H(t) \begin{bmatrix} \delta_x \\  \delta_u \end{bmatrix}\right|^2ds\ge \left|\int_0^1 H(t) \begin{bmatrix} \delta_x \\  \delta_u \end{bmatrix}ds\right|^2
\]
and if the right-hand side is computed along any path, it is an upper bound for
\[
\left| H(t) \begin{bmatrix} (x(t)-x^\star(t)) \\  (u(t)-u^\star(t)) \end{bmatrix}\right|^2.
\]
with equality if a straight-line path is used, i.e. a geodesic with respect to $|\delta_x|$.

$\Box$

As approximate optimal control, one could consider maximizing, e.g. the smallest eigenvalue of $ W(x)$ or the trace of $W^{-1}(x)$, both of which are concave in $W(x)$, to minimize ``worst case'' or ``expected'' value of 
\[
\int_0^1\delta_x' W^{-1}(x,t)\delta_x ds.
\]

The above construction can be easily extended to cost functions with a cross term $(u(t)-u^\star(t))'S(t)(x(t)-x^\star(t))$ as long as the cost matrix remains positive-definite.

It is well known that for linear systems, stabilizability implies solvability of the LQR problem if $Q>0, R>0$ \cite{hespanha2009linear}. In fact, for nonlinear systems a similar claim can be made in terms of guaranteed cost controllability.

\begin{theorem}
Suppose a system of the form \eqref{eq:sys} is universally stabilizable, as verified by the conditions of Theorem \ref{Wthm}, then there exists a solution to the quadratic guaranteed cost control problem in Theorem \ref{QGCCthm}.
\end{theorem}

\textbf{Proof:} Since both of the conditions are pointwise, it suffices to show at each state and time $(x,t)$ that the existence of a $W(x,t)$ and $\rho(x,t)$ satisfying \eqref{eq:Wcond}and 
\begin{equation}\label{eq:hypoth}
-\dot W+WA'+AW-\rho BB'<0
\end{equation}
implies the existence of a (possibly different) $\bar W(x,t)$ satisfying 
\begin{equation}\label{eq:conseq}
-\dot{\bar W}+\bar WA'+A\bar W-BR^{-1}B'+\bar W Q \bar W\le 0,
\end{equation}
for a fixed $Q(t)>0, R(t)>0$. Note that we have taken the Schur complement of \eqref{eq:Wcondgcc}. 

Let us first consider the condition
\[
\delta_1'\left(-\dot{\bar W}+\bar WA'+A\bar W-BR^{-1}B'+\bar W Q \bar W\right)\delta_1\le 0,
\]
when $\delta_1$ is in the nullspace of $B'$, i.e. $B'\delta_1=0$. Clearly this is equivalent to the condition
\begin{equation}\label{eq:p1}
\delta_1'\left(-\dot{\bar W}+\bar WA'+A\bar W+\bar W Q \bar W\right)\delta_1\le 0.
\end{equation}
Now, the conditions of Theorem \ref{Wthm} imply that for any $\delta_1$ in the null-space of $B'$, we have 
\begin{equation}\label{eq:p2}
\delta_1'\left(-\dot{ W}+WA'+A W\right)\delta_1<0.
\end{equation}
Take $\bar W(x,t) = \mu(x,t) W(x,t)$ for some scaling factor $\mu(x,t)>0$, then \eqref{eq:p1} becomes
\[
\delta_1'\left(\mu(-\dot{ W}+WA'+A W)+\mu^2WQW\right)\delta_1\le 0.
\]
From \eqref{eq:p2} it is clear that taking $\mu$ sufficiently small this condition is satisfied.

On the other hand, suppose $\delta_2$ is not in the null-space of $B$, then as $\mu\rightarrow 0$ 
\begin{align}
&\delta_2'\left(\mu(-\dot{ W}+WA'+A W)-BR^{-1}B' + \mu^2WQW\right)\delta_2\notag\\
&\rightarrow -\delta_2'BR^{-1}B'\delta_2<0\notag
\end{align}
so again, sufficiently small $\mu$ will satisfy the conditions of Theorem \ref{QGCCthm}. Since it is sufficient to verify the contraction condition for $\delta:|\delta|=1$, a compact set, this implies the existence of a sufficiently small $\mu(x,t)$ for all $\delta$.

$\Box$

\begin{remark} Taking $\bar W(x,t) = \mu W(x,t)$ very small in the proof above leads to a upper bound on the cost given by the path integral of $\frac{1}{\mu}\delta'W^{-1}(x,t)\delta$. Therefore the smaller $\mu$ is taken, the more conservative the upper bound on cost.
\end{remark}

\section{Computational Approaches} \label{sec:comp}

The methods we have described naturally break the control problem down into two stages: firstly, a pointwise LMI must be solved for $W(x)$, giving the metric $\delta_x'M(x)\delta_x$, guaranteeing a particular form of stabilizability; secondly, the on-line computation involves computing a path integral along a line or geodesic. A full discussion is not possible due to space restrictions, but in this section we briefly discuss some applicable techniques.

\subsection{Finding a Control Contraction Metric}

The main calculation is to find a matrix function $W(x,t)$ and a scalar function $\rho(x,t)$ satisfying \eqref{eq:Wcond} and \eqref{eq:Wdotcond} for all $x, u, t$, or the similar conditions in Sections \ref{sec:exp} and \ref{sec:gcc}. These conditions are convex, but the decision variables $W$ and $\rho$ are infinite dimensional and there are infinitely many LMI constraints \eqref{eq:Wcond} and \eqref{eq:Wdotcond}, due to the dependence on $x$ and $t$. To be practically computable, we must reduce these to a finite-dimensional optimization problem.

If $f(x,t)=f(x)$ is a vector of polynomials in $x$, one can search for polynomial $W(x)$ and $\rho(x)$ satisfying the conditions of the theorem using sum-of-squares programming \citep{parrilo2003semidefinite}. In that framework, a matrix inequality $G(x)\ge 0$ for all $x$, where $G(x)$ is a $q\times q$ symmetric matrix, is represented by introducing an auxiliary variable $y\in\mathbb R^q$ and search for a representation
$
y'G(x)y = \sum_{i} g_i(x,y)^2
$
for some polynomials $g_i(x,y)$. When $G(x)$ is linearly parameterized by unknowns, this can be represented as a finite-dimensional semidefinite program. This gives a conservative test, since not all non-negative polynomials are sums of squares.

Alternatively, over a compact set in state space one could specify a finite set of basis functions for $W$ and $\rho$, and check the inequalities  \eqref{eq:Wcond} and \eqref{eq:Wdotcond} at a finite grid of points. This again leads to a finite-dimensional set of linear matrix inequalities. A similar method was investigated by \cite{johansen2000computation} for computation of Lyapunov functions.

\subsection{Feedback Control via Geodesic Computation}

In the case that $D(x,\delta_x,t)=|\delta_x|$ is used to compute the distance function, construction of the feedback control is relatively straightforward: shortest paths are straight lines, and the control is computed by a path integral along a line. In the simplest case that a $M(t)$ and $\rho(t)$ can be found independent of $x$ and satisfying the conditions, then the feedback control is simply the linear gain \[u(t) =u^\star(t) -\frac{1}{2}\rho(t)B(t)'M(t)(x(t)-x^\star(t)).\]

It is to be expected that better performance can be achieved if $D(x,\delta_x,t)=\sqrt{V(x,\delta_x,t)}$ is used instead. In that case, the presented approach does not remove the need for on-line computation but it does potentially reduce its complexity, as compared to model predictive control. In particular, we have reduced the problem from an {\em optimal control} problem to a {\em shortest path} problem, with respect to a Riemannian distance function. The shortest path problem is easier to solve because of the special structure of the distance metric and the fact that the path need not satisfy any differential equations. Furthermore, it is computed only in the state space, rather than in the state/control space, as with MPC when using collocation or multiple shooting.

Algorithms for computing geodesics have been developed in several fields, including computational physics, computer graphics, and robot motion planning
\cite{memoli2001fast},
\cite{boykov2003computing}, \cite{lavalle2006planning}. The approach described in this subsection may also be applicable to the controller construction in \citep{lohmiller1998contraction}.

\subsection{Open-loop Control via Path Images}\label{sec:forwardimage}

An alternative method giving an open-loop control signal that drives $x(t)\rightarrow x^\star(t)$ was suggested in \citep{lohmiller1998contraction}. This method also uses path integrals of a differential feedback gain $\delta_u = K\delta_x$, although $K$ was constructed by a different method related to the controllability conditions of \cite{isidori1995nonlinear}. In particular, it was suggested to compute a path from $x^\star(0)$ to $x(0)$, and integrate $\delta_u = K\delta_x$ along this path for $u(0)$. For times $t>0$, the forward image of each point on the path under the flow of the system is simulated, generating a family of paths $\gamma(s,t)$, and controls $u(s,t)$.

This method can also be applied for the systems we study with $K=-\frac{1}{2}\rho B'M$. One advantage is that it does not require recomputing a geodesic at each time to compute the control signal. A disadvantage is that it does not give a feedback law, and could be non-robust. E.g., if the true system behaviour does not match the simulated forward image of $x(0)$ then there are no guarantees of stability.

One can also imagine a hybrid of this approach and the previous one, in which forward images of the initial path are computed, but with some corrective term based on gradient descent to converge towards geodesics. In this case, the resulting feedback system would be a dynamic state feedback with an infinite-dimensional controller state.

\subsection{Two-stage Computation of Explicit Feedback}

If $W(x,t)$ has been found using the above-described convex optimization procedure, then one can fix $M(x,t) = W(x,t)^{-1}$ which implies the existence of matrix $K(x)$ such that
\begin{equation}\label{eq:explcontrl}
 \delta_x'\left(\dot M(x,t) +2M(x,t)\left[A(x,t)+B(t)K(x,t)\right]\right)\delta_x<0
\end{equation}
for all $\delta_x$, i.e., the existence of a differential state feedback $\delta_u = K(x)\delta_x$. However, in order to compute an explicit control law $u=k(x,t)$, this differential relation must be integrable.

Let $K_i(x)$ refer to the $i^{th}$ row of the matrix $K(x)$, then the additional convex constraint
\begin{equation}\label{eq:Poincare}
\pder{x}K_i(x) = \left(\pder{x}K_i(x) \right)', i = 1, 2, ... m
\end{equation}
ensures the existence of a function $k(x)$ satisfying $\pder[k]{x}=K(x)$ by the Poincar\'e lemma. For example, consider a linearly parameterized class of controllers where the $i^{th}$ control element is given by
\[
u_i = k_i(x) := \sum_{j = 1}^p \kappa_{i,j} \phi_j(x)
\]
where $\phi_j(x)$ are smooth basis functions (e.g. a monomial basis for polynomial feedback). In this case,
\[
K_i(x) := \sum_{j = 1}^p \kappa_{i,j} \pder{x}\phi_j(x)
\]
and \eqref{eq:Poincare} imposes a linear constraint on the coefficients, as long as the set of basis functions is closed under differentiation with respect to $x$. This can also be considered as a curl condition on a vector field, so recent work on curl-free wavelets by \cite{deriaz2009orthogonal} may be applicable when choosing basis functions for the control gain.

For stabilizing particular equilibrium $x^\star$, one can add the additional linear constraint
\[
0 = f(x^\star,t)+B(t)k(x^\star,t),
\]
or, e.g. a periodic solution could be stabilized with family of constraints over $[0, T]$
 \[
 \dot x^\star= f(x^\star,t)+B(x^\star,t)k(x^\star,t).
 \]
It is important to note that the conditions of the main theorems imply the existence of a $K(x,t)$ satisfying \eqref{eq:explcontrl} but not necessarily satisfying the additional constraint \eqref{eq:Poincare}. To our knowledge, whether or not this stronger implication is true is an open question.

\section{Discussion}

\subsection{Necessity for the Case of Linear Systems}

In general, the conditions we provide will be sufficient but not necessary. An interesting question is whether there is a class of nonlinear systems for which they are necessary, i.e. they completely describe stabilizability. As a basic requirement, one would expect this class to include all stabilizable linear time-invariant (LTI) systems:
\[
\dot x = Ax+Bu
\]
and indeed this is the case. Our condition for LTI systems with a matrix $W>0$ independent of $x$ is:
\[
AW + WA' - BB' <0
\]
where $\rho$, also constant,has been absorbed into $W$. This is well-known Lyapunov condition for stabilizability of linear systems -- see, e.g., \citep[Sec 14.5]{hespanha2009linear}.

The condition we give on computing a guaranteed cost control similarly reduces to the Riccati inequality:
\[
AW + WA' -BR^{-1}B' + WQW\le 0
\]
which is known to have a solution if and only if the algebraic Riccati equation:
\[
AW + WA' -BR^{-1}B' + WQW = 0
\]
has a solution, given by the solution of the Riccati inequality with maximal trace.  By multiplying the latter on either side by $M = W^{-1}$, this is the exactly the Riccati inequality associated with the linear quadratic regulator problem, with the minor strengthening that $Q$ must be invertible and hence positive-definite. In linear control design, changes of variables such as that from $M$ to $W$ are frequently used to construct linear matrix inequality conditions \citep{boyd1994linear}.

\subsection{Connection to Control Lyapunov Functions}

The distance function $d(x^\star(t), x(t),t)$ is a control Lyapunov function for a solution $x^\star$, thus the conditions we give imply the existence of a control Lyapunov function for {\em every} solution of the system.

Control Lyapunov functions for systems of the form \eqref{eq:sys} can be characterized by the condition:
\[
\pder{x}V(x,t)B(t) = 0 \Rightarrow \pder{x}V(x,t)f(x,t)+ \pder{t}V(x,t) <0
\]
for all $x$ and $t$. In general, it is computationally challenging to search for a $V(x,t)$ satisfying such a condition implication. If we instead consider a similar condition for the {\em differential} dynamics \eqref{eq:diffdyn} and a function $\delta_x'M(x,t)\delta_x$:
\[
MB\delta_x = 0 \Rightarrow \delta_x'(\dot M + A'M+MA)\delta < 0
\]

Then by Finsler's theorem \citep{uhlig1979recurring}, or a version of the S-Procedure losslessness theorem \citep{yakubovich1971s}, this is equivalent to the existence of a non-negative function $\rho(x, t)$ for which the following is true:
\[
\dot M + A'M + MA -\rho MBB'M< 0
\]
which, by pre and post-multiplication by $W$ is our condition \eqref{eq:Wdotcond}. Therefore our condition guarantees both the existence of a control Lyapunov function in the traditional sense, and a ``differential'' control Lyapunov function -- i.e. a control contraction metric.

%\subsection{A Dual Perspective}
%
%The LMI construction can be interpreted in terms of the duality of linear control and filtering problems, and is closely analogous to the Lyapunov duality of density functions \cite{rantzer2001dual} and occupation measures \cite{lasserre2008nonlinear}. The latter, when considered in terms of moment sequences, also result in LMI conditions. The interpretation for the condition on $W(x)$ is that if a large number of initial conditions are sampled nearby a particular solution, then the variance of this distribution decreases over time, leading to an accumulation of probability density on the trajectory itself. If this holds for all trajectories, then all trajectories must converge.

\subsection{On Transverse Contraction}

The results in this paper also suggest a new interpretation of limit cycle stability. In \citep{manchester2013transverse} the authors introduced the concept of transverse contraction, a generalisation of prior work by \cite{Borg}, \cite{hartman1962global}, \cite{leonov1996frequency}, to study the existence of limit cycles in autonomous systems:
\[
\dot x = f(x).
\]
Transverse contraction refers to the existence of a metric function $V(x, \delta_x)$ such that 
\begin{equation}
\frac{\partial V(x,\delta)}{\partial x}f(x) + \frac{\partial V(x,\delta)}{\partial \delta}\pder[f(x)]{x}\delta\le -\lambda V(x,\delta),\label{eq:contr1}
\end{equation}
for all $\delta\ne 0$ such that $\frac{\partial V}{\partial \delta}f(x)=0$. In \citep{manchester2013transverse} it was shown that, for a Riemannian metric $\sqrt{\delta'M(x)\delta}$, this is equivalent to a pointwise LMI condition:
\begin{align}
\dot W(x) &\ge W(x)A(x)' + A(x)W(x) -\rho(x)Q(x)'+2\lambda W(x),
\end{align}
where $A(x) = \pder[f]{x}$, $Q(x)=f(x)f(x)'$ and $W(x)>0$. It was shown that this implies stability under time reparametrisation -- a.k.a. Zhukovsky stability -- which is known to imply existence of a stable limit cycle.

In the context of this paper, that condition is equivalent to stabilizability of the system:
\[
\dot x = f(x) + f(x)u.
\]
The interpretation of this fact is that a control input acting ``along'' the flow of the system $f(x)$ can speed up or slow down the trajectory, but not change its phase portrait.

\section{Extension to more General Nonlinear Systems} \label{sec:gen}

In this section we briefly discuss extension to more general nonlinear systems
\[
\dot x(t) = f(x(t),u(t), t)
\]
with differential dynamics
\[
\dot \delta_x(t) = A(x,u,t)\delta_x(t) + B(x,u,t)\delta_u(t).
\]
The significant change from the previous sections is that the $A$ and $B$ matrices now depend on $u$.

In this case, we can also consider a matrix $W(x,u,t)$ and constant $\rho(x,u,t)$. Then, suppose conditions \eqref{eq:Wcond} and \eqref{eq:Wdotcond} hold for all $x, u, t$. The path integral equation for the control signal $u(s)$ now depends on $u(s)$, and can be constructed if a solution exists to the differential equation:
\[
\frac{d}{ds}u = -\frac{1}{2} \rho(\gamma,u)B(\gamma)'M(\gamma, u)\pder[\gamma]{s}(s)
\]
with boundary condition $u(0) = u^\star(t)$. If the dependency of these terms on $u$ can be well controlled, then it may be possible to guarantee existence by, e.g., the Bellman-Gr\"onwall lemma. We defer detailed discussion of this issue for a later publication.

\section{Application Example}

In this section we show an example class of systems for which the above analysis can be applied: compositions of a stabilizable linear system with a contracting system.

Consider the following hierarchical composition of systems:
\[
\dot z = \bar Az+\bar Bu, \ \ \dot y = f_1(y,z),
\]
for which $f_1$ is partially contracting with respect to $y$, and the pair $\bar A, \bar B$ is stabilizable. That is, there exists $M_1(y)\ge \alpha I$ and $P$ such that
\begin{align}
\dot M_1(y) + \pder[f_1]{y}'M_1(y)+M_1(y)\pder[f_1]{y}<&0 \label{eq:seriescontr}\\
\bar AP+P\bar A'-\bar B\bar B'<&0.\label{eq:seriesctrb}
\end{align}
Let $W_1(y) = M_1^{-1}(y)$, then we also have from \eqref{eq:seriescontr}
\begin{equation}\label{eq:seriescontr2}
\Pi(y):=-\dot W_1(y) + W_1(y)\pder[f_1]{y}'+\pder[f_1]{y}W_1(y)<0
\end{equation}
with $W_1(y) \le \frac{1}{\alpha} I$.

The combined system with state $x:=[y', z']'$ is:
\[
\dot x =\begin{bmatrix}\dot y\\ \dot z\end{bmatrix} =\begin{bmatrix}f_1(y,z)\\ \bar Az\end{bmatrix} + \begin{bmatrix}0\\\bar B\end{bmatrix} u=:f(x)+Bu
\]
and the differential dynamics are \eqref{eq:diffdyn} with
\[
A(y) := \pder[f]{x} = \begin{bmatrix} \pder[f]{y} & \pder[f]{z} \\ 0 & \bar A \end{bmatrix}
\]
Let us consider a class of metrics parametrised by $\beta>0$:
\[
W(y)=\begin{bmatrix}\beta(y) M_1(y)^{-1}&0\\0&P\end{bmatrix}
\]
Noting that
\begin{align}
\dot W =
\begin{bmatrix}
\dot W_1&0\\0&0
\end{bmatrix},
BB'=
\begin{bmatrix}
0&0\\0&\bar B \bar B'
\end{bmatrix}
\end{align}
we see that \eqref{eq:Wdotcond} reduces to
\begin{equation}
\begin{bmatrix} \beta(y)\Pi(y)& \beta(y) W_y\pder[f]{z}\\ \beta(y) \pder[f]{z}'W_y& \bar A'P+P\bar A-\bar B\bar B'\end{bmatrix}<0.
\end{equation}
Now, since the top-left block is negative definite by assumption \eqref{eq:seriescontr2}, negativity of the entire matrix is equivalent, by the Schur complement, to the following
\begin{align}
\bar A'P+P\bar A-\bar B\bar B'-
\beta(y) W_y\pder[f]{z}\Pi(y)\pder[f]{z}'W_y<0.
\end{align}
Now, since $\bar A'P+P\bar A-\bar B\bar B'<0$ by assumption \eqref{eq:seriesctrb}, this can clearly be satisfied for sufficiently small $\beta(y)>0$. Thus we have another class of systems for which the conditions of Theorem \ref{Wthm} are necessary and sufficient.

This situation is clearly simpler than Lyapunov-based control design: if a nonlinear system is Lyapunov stable at the origin, and driven by a stabilizable linear system, it is not  so clear that the entire system is stabilizable. This is the fact that led to the development of backstepping \citep{krstic1995nonlinear}.

\section{Conclusions}

In this paper we have introduced the concept of {\em universal stabilizability}: the condition that {\em every} solution of a system is globally stabilizable. We have proven sufficient conditions in terms of the existence of a {\em control contraction metric}.

Unlike control Lyapunov functions, the set of control contraction metrics for a given system can be parametrised as a convex set -- defined by pointwise linear matrix inequality constraints -- and thus amenable to search via convex optimization methods. 

The conditions we give are necessary and sufficient for linear systems and certain classes of interconnected nonlinear systems. Straightforward extensions allow one to construct convex upper bounds for a nonlinear quadratic regulator problem.

%
%Condition for stabilizability
%
%Extension to exp stab, LQ sub-optimality.
%
%Construction of explicit control laws
%

\bibliographystyle{alpha}
\bibliography{elib}

\end{document}